 \documentclass[final,1p,times]{elsarticle}

\usepackage{amsmath,amsthm,amssymb,enumitem,hyperref}
\usepackage{amsfonts}
\usepackage{textcomp}
\usepackage{xcolor}
\usepackage[utf8]{inputenc}

\newtheorem{theorem}{Theorem}[section]

\newtheorem{prop}[theorem]{Proposition}

\newcommand{\bracenom}{\genfrac{\lbrace}{\rbrace}{0pt}{}}

\usepackage{ulem}

\begin{document}

\begin{frontmatter}



\title{Sign-consistent estimation in a sparse Poisson model}


\author{M. Gomtsyan, C. L\'evy-Leduc, S. Ouadah, L. Sansonnet}

\affiliation{organization={Université Paris-Saclay, AgroParisTech, INRAE, UMR MIA Paris-Saclay},
            addressline={22 place de l'agronomie}, 
            city={Palaiseau},
            postcode={91120}, 
            country={France}}

\begin{abstract}
In this work, we consider an estimation method in sparse Poisson models inspired by \cite{friedman_hastie_tibshirani:2010} and provide novel sign consistency results under mild conditions.  
\end{abstract}




\end{frontmatter}


\section{Introduction}

Discrete-valued data arise in diverse applied scientific areas, ranging from finance to molecular biology and epidemiology. 
For example, as discussed in~\cite{wu2017diversity}, in molecular biology, non-coding genes are emerging as potential key regulators of the expression of protein-coding genes. Yet, among numerous non-coding genes, only a few are likely to be involved for explaining the expression of the coding genes. Consequently, variable selection will help to identify the relevant non-coding genes by obtaining sparse estimators, meaning that most of them are zero. A  popular approach in statistics for performing variable selection is the Lasso proposed by \cite{tibshirani:1996}. 
Besides, \cite{zhao:2006} showed that Lasso has theoretical guarantees, under some mild conditions. More particularly, Lasso is model selection consistent, meaning that Lasso chooses the true model. However, the consistency results are established in a Gaussian setting, which may not hold for discrete-valued data. 

In this work we consider the following sparse Poisson model. Let $Y_1, \ldots, Y_n$ be independent random variables such that for all $i$, 
\begin{equation}
Y_i\sim\texttt{Poisson}(\lambda^{\star}_i) \quad \text{with} \quad \lambda^{\star}_i = \exp(x_i\pmb{\beta}^{\star}), 
\label{eq_model}
\end{equation}
where $x_i$ is the $i$th row of a $n \times p$ design matrix $\mathbf{X}$ and $\pmb{\beta}^{\star}$ is a sparse vector of regression coefficients in $\mathbb{R}^p$. The non-null coefficients correspond to the predictors that are relevant to explain the response. In the following we will consider an estimation method inspired by \cite{friedman_hastie_tibshirani:2010} and provide a novel sign consistency result.

The paper is organised as follows. Firstly, in Section 2 we present the statistical approach for estimating $\pmb{\beta}^{\star}$. Next, in Section 3 we establish its sign-consistency. The detailed proof is available in the Appendix.

\section{Statistical approach}

To estimate $\pmb{\beta}^{\star}$ defined in Model (\ref{eq_model}), we shall use the approximation proposed by \cite{friedman_hastie_tibshirani:2010} which consists in maximizing with respect to $\pmb{\beta}$ the second order Taylor approximation of the log-likehood $l$ at the current estimate $\pmb{\Tilde{\beta}}$, namely:
\begin{equation*}
  l(\pmb{\Tilde{\beta}}) +  \sum_{i=1}^n \sum_{k=1}^p{x_i}_k (Y_i - \Tilde{\lambda}_i) (\beta_k - \Tilde{\beta}_k)
  -\frac{1}{2}\sum_{i=1}^n  \sum_{1\leq k,\ell\leq p}  \Tilde{\lambda}_i  {x_i}_k (\beta_k - \Tilde{\beta}_k) {x_i}_\ell (\beta_\ell - \Tilde{\beta}_\ell),
\end{equation*}   
where $\Tilde{\lambda}_i = \exp(x_i  \pmb{\Tilde{\beta}})$, 
which means maximizing
\begin{equation*}
  \sum_{i=1}^n(Y_i - \Tilde{\lambda}_i) x_i(\pmb{\beta} - \pmb{\Tilde{\beta}})-\frac{1}{2}\sum_{i=1}^n \Tilde{\lambda}_i
  \left({x_i}(\pmb{\beta} - \pmb{\Tilde{\beta}}) \right)^2
=  \sum_{i=1}^n\frac{(Y_i - \Tilde{\lambda}_i)}{\sqrt{\Tilde{\lambda}_i}} \left(\sqrt{\Tilde{\lambda}_i}x_i(\pmb{\beta}- \pmb{\Tilde{\beta}})\right) - \frac{1}{2} \sum_{i=1}^n  \left(\sqrt{\Tilde{\lambda}_i}{x_i}\left(\pmb{\beta} - \pmb{\Tilde{\beta}}\right) \right)^2.
\end{equation*}
This boils down to minimizing
$$
\sum_{i=1}^n \left(\sqrt{\Tilde{\lambda}_i}{x_i}\left(\pmb{\beta} - \pmb{\Tilde{\beta}}\right) -\frac{Y_i - \Tilde{\lambda}_i}{\sqrt{\Tilde{\lambda}_i}}\right)^2.
$$
Minimizing this criterion can be viewed as
the minimization with respect to $\pmb{\beta}$ of the following least-squares criterion:
$\|\mathcal{Y} - \mathcal{X} \pmb{\beta} \|_2^2$ where $\|u\|_2^2=\sum_{i=1}^n u_i^2$ for a vector $u=(u_1,\dots,u_n)$ in $\mathbb{R}^n$,
\begin{equation}\label{eq:Y_X}
   \mathcal{Y}=\mathcal{X}\pmb{\Tilde{\beta}}+\boldsymbol{\Tilde{\Lambda}}^{-1/2}(\textbf{Y}-\boldsymbol{\tilde{\lambda}})
   \quad \textrm{ and }\quad
   \mathcal{X}=\boldsymbol{\Tilde{\Lambda}}^{1/2}\textbf{X},
\end{equation}
 $\boldsymbol{\Tilde{\Lambda}}$ denoting the diagonal matrix having the
 $\Tilde{\lambda}_i$'s as diagonal elements, $\boldsymbol{\tilde{\lambda}}$ being a column vector having the $\Tilde{\lambda}_i$'s as components
 and $\textbf{Y}$ denoting a column vector having the $Y_i$'s as components.
   
Thus, in order to obtain a sparse estimation of $\pmb{\beta}^\star$, we will focus on finding $\pmb{\hat{\beta}}(\alpha)$ defined for $\alpha>0$ by:
\begin{equation}\label{eq:def_beta_hat}
\pmb{\hat{\beta}}(\alpha) = \underset{\pmb{\beta}\in\mathbb{R}^p}{\arg\min} \Bigg\{ \|\mathcal{Y} - \mathcal{X} \pmb{\beta}\|_2^2 + \alpha \|\pmb{\beta}\|_1 \Bigg\},
\end{equation}
where $\|v\|_1=\sum_{k=1}^p|v_k|$ for a vector $v=(v_1,\dots,v_p)$ in $\mathbb{R}^p$.

We shall establish the sign consistency of $\pmb{\hat{\beta}}$ in Theorem \eqref{theo1} of the following section.

\section{Sign consistency}
%
Let 
\begin{equation}\label{eq:def_C_W}
\textbf{C} = \frac{\mathcal{X}^T \mathcal{X}}{n} \quad \textrm{and}\quad 
    \textbf{W} = \frac{\mathcal{X}^T \pmb{\Tilde{\varepsilon}}}{n} ,
\end{equation}
where $A^T$ denotes the transpose of the matrix $A$,
\begin{equation}
\boldsymbol{\Tilde{\varepsilon}} = (\Tilde{\varepsilon}_1, \dots, \Tilde{\varepsilon}_n)^T \quad \text{with} \quad \Tilde{\varepsilon}_k = \frac{Y_k - \Tilde{\lambda}_k}{\sqrt{\Tilde{\lambda}_k}}  \quad \text{for all } \quad 1 \leq k \leq n.
\label{eq_eps}
\end{equation}
Without loss of generality, suppose that $\pmb{\beta}^{\star} = (\beta^{\star}_1, \dots \beta^{\star}_q, \beta^{\star}_{q+1}, \dots \beta^{\star}_p)^T$, where $\beta^{\star}_j \neq 0$ when $1 \leq j \leq q$ and $\beta^{\star}_j =0$ when $q+1 \leq j \leq p$ and denote
\begin{equation}\label{eq:beta_star}
    \pmb{\beta}_1^{\star} = (\beta^{\star}_1, \dots, \beta^{\star}_q)^T \quad \text{and} \quad  \pmb{\beta}_2^{\star} = (\beta^{\star}_{q+1}, \dots, \beta^{\star}_p)^T.
  \end{equation}
Then,
\begin{equation}\label{eq:def_C_W_1_2}
\textbf{C} = 
\begin{pmatrix}
\mathcal{X}_1^T \mathcal{X}_1 / n & \mathcal{X}_1^T \mathcal{X}_2 / n \\
\mathcal{X}_2^T \mathcal{X}_1 / n & \mathcal{X}_2^T \mathcal{X}_2 / n \\
\end{pmatrix} = 
\begin{pmatrix}
C_{11} & C_{12} \\
C_{21} & C_{22} \\
\end{pmatrix}
\quad 
\textrm{ and }
\quad
\textbf{W} = 
\begin{pmatrix}
W_1 \\
W_2 \\
\end{pmatrix}.
\end{equation}
\begin{theorem}\label{theo1}
 Assume that $Y_1, \dots, Y_n$ are independent random variables such that for all $i$, $Y_i\sim\texttt{Poisson} (\lambda^{\star}_i)$ with $\lambda^{\star}_i = \exp(x_i\pmb{\beta}^{\star})$, where
$x_i$ is the $i$th row of a design matrix $\mathbf{X}$ and $\pmb{\beta}^{\star}$ is defined in (\ref{eq:beta_star}).
  Assume also that there exist positive constants $M_1$, $M_2$, $M_3$, $M_4$, $M_5$, $M_6$, $M_7$ and $c_1$  such that $0 < c_1 \leq 1$, that
  $p$ does not depend on $n$ and that the following assumptions hold. 
\begin{enumerate}[label=(T\arabic*)]
\item \label{th_1} For all $k$ in $\{1,\dots,n\}$, $\|x_k \|_2 \leq M_1$ and for all $\ell$ in $\{1,\dots,p\}$, $\|x^{(\ell)}\|_2\leq M_7$, where
  $x^{(\ell)}$ denotes the $\ell$th column of $\mathbf{X}$.
\item \label{th_5} $\pmb{\Tilde{\beta}}$ is a preliminary estimator of $\pmb{\beta}^\star$ such that
    $|\beta^\star_i - \Tilde{\beta}_i| =O_P(1/n)$, as $n$ tends to infinity, for all $i=1, \dots, p$.
\item \label{th_2} With a probability tending to 1 as $n$ tends to infinity, $\lambda_{\min}(C_{11}) \geq M_2$, where $\lambda_{\min}(A)$ denotes the smallest eigenvalue of the matrix $A$ and $C_{11}$ is defined in (\ref{eq:def_C_W_1_2}).
  \item \label{th_3} With a probability tending to 1 as $n$ tends to infinity, $\lambda_{\max}(C_{12}) \leq M_3$, $\lambda_{\max}(C_{21}) \leq M_4$, and $\lambda_{\max}(C_{22}) \leq M_5$, where $\lambda_{\max}(A)$
    denotes the largest eigenvalue of the matrix $A$.
  \item \label{th_4} $n^{\frac{1-c_1}{2}} \min_{1 \leq i \leq q } |\beta_i^{\star}| \geq M_6$.
\end{enumerate}
Let us also suppose that the following condition called strong irrepresentable condition holds: there exists  $\tau > 2/3$ such that
\begin{equation}\label{eq:IC}
\left| C_{21} C_{11}^{-1}\text{sign}(\pmb{\beta}^\star_1) \right| \leq 1- \tau,
\end{equation}
with a probability tending to 1 when $n$ tends to infinity, where the inequality has to be understood component by component. 
Then, for all $\alpha=\alpha_n$ such that
$
\alpha_n=O\left(n^\frac{c_2+1}{2} \right), \quad \text{where} \quad 0 < c_2 < c_1 \leq 1$, $\pmb{\hat{\beta}}$ defined in (\ref{eq:def_beta_hat}) satisfies
\begin{equation*}
\mathbb{P} \left( \text{sign}(\pmb{\hat{\beta}}) = \text{sign}(\pmb{\beta}^\star)\right) \rightarrow 1, \quad \text{when} \quad n \rightarrow \infty.
\end{equation*}
\end{theorem}

The proof of Theorem \ref{theo1} relies on the following proposition.


\begin{prop}\label{prop1}
  Under the assumptions of Theorem \ref{theo1}, let
\begin{equation}\label{eq:R1_2}
R_1 = (C(\pmb{\beta}^{\star} - \pmb{\Tilde{\beta}}))_1 \quad \text{and} \quad R_2 = (C(\pmb{\beta}^{\star} - \pmb{\Tilde{\beta}}))_2.
\end{equation}
Then,
\begin{equation*}
    \mathbb{P}\Big(\text{sign}(\pmb{\hat{\beta}})=\text{sign}(\pmb{\beta}^{\star})\Big) \geq \mathbb{P}\Big(A_n \cap B_n \Big),
\end{equation*}
where
\begin{equation*}
    A_n = \Big\{ |C_{11}^{-1} W_1| < |\pmb{\beta}_1^{\star}| - \frac{\alpha}{2n} |C_{11}^{-1} \text{sign}(\pmb{\beta}_1^{\star}) | - |C_{11}^{-1} R_1| \Big\}
\end{equation*}
and
\begin{equation*}
    B_n = \left\{ | C_{21} C_{11}^{-1} W_1 - W_2 | \leq \frac{\alpha}{2n} \left( 1 - | C_{21}C_{11}^{-1} \text{sign}(\pmb{\beta}_1^{\star}) | \right) - | C_{21}C_{11}^{-1} R_1 - R_2 |  \right\}.
\end{equation*}
\end{prop}


The proof of Proposition \ref{prop1} is in \ref{sec:proof_prop1}.

Proving Theorem \ref{theo1} consists in showing that $\mathbb{P}(A_n^c)$ and $\mathbb{P}(B_n^c)$ go to zero as $n$ tends to infinity
where $S^c$ denotes the complementary set of the set $S$. The proof is given in \ref{sec:proof_theo}.

\appendix 
\section{Proofs}
\subsection{Proof of Proposition \ref{prop1}}\label{sec:proof_prop1}

First observe that we get from (\ref{eq:Y_X}) and (\ref{eq_eps}) that
\begin{equation}\label{eq:def_Y}
\mathcal{Y}=\mathcal{X}\pmb{\Tilde{\beta}}+\boldsymbol{\Tilde{\varepsilon}}.
\end{equation}
Let us denote $\pmb{\hat{\beta}}$ the estimator $\pmb{\hat{\beta}}(\alpha)$ defined in (\ref{eq:def_beta_hat}), it satisfies the following Karush–Kuhn–Tucker conditions described in \cite[Section 4.2.2]{giraud:2021}:
\begin{center}
\begin{tabular}{cl}
$\displaystyle \big( \mathcal{X}^T(\mathcal{Y} - \mathcal{X} \pmb{\hat{\beta}}) \big)_i = \frac{\alpha}{2} \text{sign}(\hat{\beta}_i)$, & if $\hat{\beta}_i \neq 0$, \\
$\displaystyle \big| \big( \mathcal{X}^T (\mathcal{Y} - \mathcal{X} \pmb{\hat{\beta}}) \big)_i \big| \leq \frac{\alpha}{2}$, & if $\hat{\beta}_i = 0$,
\end{tabular}
\end{center}
which can be rewritten as follows by using (\ref{eq:def_C_W}) and (\ref{eq:def_Y})
\begin{align}
&  \big(\textbf{C}(\pmb{\hat{\beta}} -  \pmb{\beta}^{\star})-\textbf{W}+\textbf{C}(\pmb{\beta}^{\star} - \pmb{\Tilde{\beta}})\big)_i = - \frac{\alpha}{2n} \text{sign}(\hat{\beta}_i), \quad \text{if } \hat{\beta}_i \neq 0,\label{KKT1}\\
&  \big| \big(\textbf{C}(\pmb{\hat{\beta}} -  \pmb{\beta}^{\star})-\textbf{W}+\textbf{C}(\pmb{\beta}^{\star} - \pmb{\Tilde{\beta}})\big)_i\big| \leq \frac{\alpha}{2n}, \quad \text{if } \hat{\beta}_i = 0.\label{KKT2}
\end{align}

If $A_n$ holds then
\begin{equation}\label{eq:AN_consequence}
        -  |\pmb{\beta}_1^{\star}| < C_{11}^{-1} W_1 - \frac{\alpha}{2n} C_{11}^{-1} \text{sign}(\pmb{\beta}_1^{\star}) - C_{11}^{-1} R_1 < |\pmb{\beta}_1^{\star}|.
\end{equation}
Let $\pmb{\check{\beta}}=(\pmb{\check{\beta}}_1^T,\textbf{0}^T)^T$ where
\begin{equation}\label{eq:check_beta}
\pmb{\check{\beta}}_1 = \pmb{\beta}_1^{\star}+C_{11}^{-1} W_1 -\frac{\alpha}{2 n} C_{11}^{-1} \text{sign}(\pmb{\beta}_1^{\star}) - C_{11}^{-1} R_1.
\end{equation}
By (\ref{eq:AN_consequence}), $|\pmb{\check{\beta}}_1-\pmb{\beta}_1^{\star}|<|\pmb{\beta}_1^{\star}|$ which implies that
$\text{sign}(\pmb{\check{\beta}}_1)= \text{sign}(\pmb{\beta}_1^{\star})$. Hence, $\pmb{\check{\beta}}$ satisfies (\ref{KKT1}).

If $B_n$ holds then
\begin{equation*}
  \Big| C_{21} \Big( C_{11}^{-1} W_1 - C_{11}^{-1} R_1 -\frac{\alpha}{2 n} C_{11}^{-1} \text{sign}(\pmb{\beta}_1^{\star}) \Big) - W_2 + R_2 \Big|
  \leq \frac{\alpha}{2 n},
\end{equation*}
which by (\ref{eq:check_beta}) corresponds to (\ref{KKT2}) for $\pmb{\check{\beta}}$ and concludes the proof of Proposition \ref{prop1}.

\subsection{Proof of Theorem \ref{theo1}}\label{sec:proof_theo}

Let us first prove that $\mathbb{P}(A_n^c)$ tends to zero as $n$ tends to infinity.
By denoting 
\begin{equation}
  \xi = (\xi_1, \dots, \xi_q)^T = C_{11}^{-1} W_1  \quad \textrm{ and }\quad
  b = (b_1, \dots, b_q)^T = C_{11}^{-1} \text{sign}(\pmb{\beta}_1^{\star}),
  \label{def:xi}
\end{equation}
we get that
\begin{align*}
  \mathbb{P}(A_n^c) &= \mathbb{P} \left( |C_{11}^{-1}W_1| + |C_{11}^{-1} R_1|+ \frac{\alpha}{2n} \left| C_{11}^{-1} \text{sign}(\pmb{\beta}_1^{\star})\right|
                     \geq |\pmb{\beta}_1^{\star}|    \right) \\
                     &   \leq \sum_{j=1}^q \mathbb{P}\left(|\xi_j| + |(C_{11}^{-1} R_1)_j|+\frac{\alpha}{2n} \big| b_j \big|\geq |\beta_j^{\star}|\right)\\
  &  \leq \sum_{j=1}^q \left\{\mathbb{P}\left(|\xi_j|\geq \frac{|\beta_j^{\star}|}{3}\right)+\mathbb{P}\left(|(C_{11}^{-1} R_1)_j|\geq \frac{|\beta_j^{\star}|}{3}\right)
    +\mathbb{P}\left(\frac{\alpha}{2n} \big| b_j \big|\geq \frac{|\beta_j^{\star}|}{3}\right)\right\}.
  \end{align*}
 By the Cauchy-Schwarz inequality, we get that for all $j$ in $\{1,\dots,q\}$
\begin{equation*}
  |b_j| \leq \sum_{j=1}^{q} |b_j| \leq\sqrt{q} \|b\|_2=\sqrt{q} \| C_{11}^{-1} \text{sign}(\pmb{\beta}_1^{\star}) \|_2 \leq q \| C_{11}^{-1}\|_2
  = q \lambda_{\text{max}}(C_{11}^{-1}).
\end{equation*}
By using \ref{th_4}, \ref{th_2}, $\alpha=O(n^{(c_2+1)/2})$ and $0<c_2<c_1\leq 1$, we obtain that
\begin{equation}\label{eq:An_c_bound}
    \mathbb{P}(A_n^c) \leq \sum_{j=1}^{q} \mathbb{P} \Bigg( |\xi_j| \geq \frac{M_6 n^{\frac{c_1-1}{2} }}{3}\Bigg) + \sum_{j=1}^{q} \mathbb{P} \Bigg(  |(C_{11}^{-1} R_1)_j| \geq \frac{M_6 n^{\frac{c_1-1}{2}}}{3}\Bigg)+o(1).
  \end{equation}

Let us first prove that the second term in the r.h.s of (\ref{eq:An_c_bound}) tends to 0 as $n$ tends to infinity.
  Observing that $R_1$ defined in (\ref{eq:R1_2}) satisfies: $R_1= C_{11}\left( \pmb{\beta}^\star_1 - \pmb{\Tilde{\beta}}_1  \right) + C_{12}\left( \pmb{\beta}^\star_2 - \pmb{\Tilde{\beta}}_2  \right)$, we get by the Cauchy-Schwarz inequality that
 \begin{equation*}
\left | (C_{11}^{-1} R_1)_j \right | \leq \sqrt{q} \left \| C_{11}^{-1} R_1\right \|_2 \leq \sqrt{q} \left \| \pmb{\beta}^\star_1 - \pmb{\Tilde{\beta}}_1  \right\|_2 + \sqrt{q} \left \| C_{11}^{-1} C_{12} \left(\pmb{\beta}^\star_2 - \pmb{\Tilde{\beta}}_2  \right) \right \|_2.
\end{equation*}
By \ref{th_5}, \ref{th_2}, \ref{th_3}, $\left | (C_{11}^{-1} R_1)_j \right | =O_P(n^{-1})$, for all $j$, where the $O_P$ does not depend on $j$,
which proves that the second term in the r.h.s of (\ref{eq:An_c_bound}) tends to 0 as $n$ tends to infinity. 

Let us now prove that the first term in the r.h.s of (\ref{eq:An_c_bound}) tends to 0 as $n$ tends to infinity. By \eqref{def:xi}, \eqref{eq:def_C_W} and \eqref{eq_eps}, denoting $\boldsymbol{\lambda^{\star}}$ the column vector of the $\lambda_i^\star$'s, $\xi$ can be rewritten as follows
\begin{align}\label{eq:xi}
\xi_j & = (C_{11}^{-1} W_1)_j =  \left(C_{11}^{-1} \frac{\mathcal{X}_1^T \pmb{\Tilde{\varepsilon}}}{n}\right)_j  = \left(\frac1n C_{11}^{-1} \mathcal{X}_1^T \boldsymbol{\Tilde{\Lambda}}^{-1/2}(\textbf{Y}-\boldsymbol{\tilde{\lambda}})\right)_j\nonumber\\
 & = \left(\frac1n C_{11}^{-1} \mathcal{X}_1^T \boldsymbol{\Tilde{\Lambda}}^{-1/2}(\textbf{Y}-\boldsymbol{\lambda^{\star}})\right)_j
  +\left(\frac1n C_{11}^{-1} \mathcal{X}_1^T \boldsymbol{\Tilde{\Lambda}}^{-1/2}(\boldsymbol{\lambda^{\star}}-\boldsymbol{\tilde{\lambda}})\right)_j.
\end{align}
For all $j$ in $\{1,\dots,q\}$, all $k$ in $\{1,\dots,n\}$, by (2.3.8) of \cite{golub:96}, we have that
\begin{align}
&\frac1n\left(C_{11}^{-1} \mathcal{X}_1^T \boldsymbol{\Tilde{\Lambda}}^{-1/2}\right)_{jk}  
\leq \frac{1}{\sqrt{n}}\left|\left(C_{11}^{-1} \frac{\mathcal{X}_1^T}{\sqrt{n}} \right)_{jk}\right| \sup_{k\in\{1,\ldots,n\}}\left(\tilde{\lambda}_k^{-1/2}\right)
\leq \frac{1}{\sqrt{n}}\left\| C_{11}^{-1} \frac{\mathcal{X}_1^T}{\sqrt{n}} \right\|_2 \sup_{k\in\{1,\ldots,n\}}\left(\tilde{\lambda}_k^{-1/2}\right)\nonumber\\ 
&= \frac{1}{\sqrt{n}}\;\rho\left(\left( C_{11}^{-1} \frac{\mathcal{X}_1^T}{\sqrt{n}}\right)^T\left(C_{11}^{-1} \frac{\mathcal{X}_1^T}{\sqrt{n}}\right)\right)^{1/2} \sup_{k\in\{1,\ldots,n\}}\left(\tilde{\lambda}_k^{-1/2}\right), \label{eq:Gjk_9}
\end{align}
where $\rho(A)$ is the spectral radius of the matrix $A$. Note that, by Theorem 1.3.22 of \cite{horn_johnson:2013},
\begin{equation}\label{eq:rho}
  \rho\left(\left(C_{11}^{-1} \frac{\mathcal{X}_1^T}{\sqrt{n}}\right)^T\left(C_{11}^{-1} \frac{\mathcal{X}_1^T}{\sqrt{n}}\right)\right) = \rho\left(\left(C_{11}^{-1}
      \frac{\mathcal{X}_1^T}{\sqrt{n}}\right)\left(C_{11}^{-1} \frac{\mathcal{X}_1^T}{\sqrt{n}}\right)^T\right)
  = \rho\left(C_{11}^{-1} \frac{\mathcal{X}_1^T \mathcal{X}_1}{n} C_{11}^{-1}\right)= \rho\left(C_{11}^{-1} \right).
\end{equation}
From \eqref{th_1} and \eqref{th_5}, we get that 

$$\tilde{\lambda}_k^{-1/2}=\frac{1}{\sqrt{\exp(x_k\pmb{\tilde{\beta}})}}=\frac{1}{\sqrt{\exp(x_k\pmb{\beta^\star})\exp(x_k(\pmb{\tilde{\beta}}-\pmb{\beta^\star}))}}={\lambda_k^\star}^{-1/2}\left(1+O_P(n^{-1})\right),$$ where ${\lambda_k^\star}^{-1/2}\leq\exp(-x_k\boldsymbol{\beta^\star}/2)\leq\exp(\|x_k\|_2
\|\boldsymbol{\beta^\star}\|_2/2)\leq\exp(M_1\|\boldsymbol{\beta^\star}\|_2/2)$, by the Cauchy-Schwarz inequality and \ref{th_1}.   Hence,
\begin{equation}\label{eq:suplambda}
\sup_{k\in\{1,\ldots,n\}}\left(\tilde{\lambda}_k^{-1/2}\right)=O_P(1), \textrm{ as } n\to\infty.
\end{equation}
By (\ref{eq:Gjk_9}), (\ref{eq:rho}) and (\ref{eq:suplambda})
$$
\frac1n\left(C_{11}^{-1} \mathcal{X}_1^T \boldsymbol{\Tilde{\Lambda}}^{-1/2}\right)_{jk}=O_P(n^{-1/2}).
$$
Moreover, by \ref{th_1} and \ref{th_5}
\begin{multline}\label{eq:lambdak-lambdakstar}
  \lambda_k^\star-\Tilde{\lambda}_k=\exp(x_k\pmb{\beta^\star})-\exp(x_k\pmb{\tilde{\beta}})=\exp(x_k\pmb{\beta^\star})\left(1-\exp(x_k(\pmb{\tilde{\beta}}-
    \pmb{\beta^\star}))\right)\\
  =-\exp(x_k\boldsymbol{\beta}^\star)\sum_{\ell\geq 1} \frac{(x_k(\boldsymbol{\Tilde{\beta}}-\boldsymbol{\beta^\star}))^\ell}{\ell !}
=O_P(1/n).
\end{multline}
Thus, the second term in the r.h.s of (\ref{eq:xi}) is $O_P(1/\sqrt{n})$. 
Note that
\begin{align}\label{eq:2eme_terme}
&\left(\frac1n C_{11}^{-1} \mathcal{X}_1^T \boldsymbol{\Tilde{\Lambda}}^{-1/2}(\textbf{Y}-\boldsymbol{\lambda^{\star}})\right)_j
=\left(\frac1n C_{11}^{-1} X_1^T (\textbf{Y}-\boldsymbol{\lambda^{\star}})\right)_j\nonumber\\
&=\left(\frac1n (C_{11}^{\star})^{-1} X_1^T (\textbf{Y}-\boldsymbol{\lambda^{\star}})\right)_j+
\left(\frac1n (C_{11}^{-1} -(C_{11}^{\star})^{-1}) X_1^T (\textbf{Y}-\boldsymbol{\lambda^{\star}})\right)_j,
\end{align}
where $C^\star={\mathcal{X}^\star}^T\mathcal{X}^\star/n$ with $\mathcal{X}^\star=(\Lambda^\star)^{1/2}\mathbf{X}$, $\Lambda^\star$ being a diagonal matrix having as diagonal entries the $\lambda_k^\star$'s.
Hence, the second term in the r.h.s of (\ref{eq:2eme_terme}) is bounded by
$$
\|C_{11}^{-1} -(C_{11}^{\star})^{-1}\|_2 \; \|X_1^T\|_2 \; \left(\frac1n\sum_{k=1}^n |Y_k-\lambda_k^\star|\right)=\|C_{11}^{-1} -(C_{11}^{\star})^{-1}\|_2 \; O_P(1),
$$
by Markov's inequality and \ref{th_1}. Since, $C_{11}^{\star}=C_{11}+(C_{11}^{\star}-C_{11})$ , we have
$$(C_{11}^{\star})^{-1}=\left(\textrm{Id}+C_{11}^{-1}(C_{11}^{\star}-C_{11})\right)^{-1}C_{11}^{-1}$$ and thus by Corollary 5.6.16 of \citep{horn_johnson:2013}, we get that
\begin{equation}\label{eq:norme2_bound_C11-1}
\|C_{11}^{-1} -(C_{11}^{\star})^{-1}\|_2=O_P(\|C_{11}^{\star}-C_{11}\|_2)=O_P(1/n^2)
\end{equation}
since by Equation (2.3.8) of \cite{golub:96}, the Cauchy-Schwarz inequality and \eqref{eq:lambdak-lambdakstar}
\begin{equation}\label{eq:norme2_bound}
\|C_{11}^{\star}-C_{11}\|_2\leq \frac{q}{n} \max_{1\leq \ell\leq q}\|x^{(\ell)}\|_2^2\sup_{1\leq k\leq n}|\lambda_k^\star-\tilde{\lambda}_k| =O_P(1/n^2).
\end{equation}
Thus, the second term in the r.h.s of (\ref{eq:2eme_terme}) is $O_P(1/n^2)$. To address the first term in the r.h.s of (\ref{eq:2eme_terme}), we shall use
the following result. 

 \begin{theorem}
 \label{bernstein}
(Bernstein's Inequality, \cite[Corollary 2.11]{boucheron_lugosi_massart:2013}) Let $X_1, \dots, X_n$ be independent real random variables. Suppose that there exist $\nu >0$ and $c>0$ such that $\sum_{k=1}^n \mathbb{E}\big[X_k^2\big] \leq \nu$, and 
\begin{equation*}
\sum_{k=1}^n \mathbb{E}\big[ |X_k|^{\ell} \big] \leq \frac{\ell!}{2} \nu c^{\ell-2}
\end{equation*}
for all integers $\ell \geq 3$. Then, for all $t>0$,
\begin{equation*}
\mathbb{P} \Bigg( \left|\sum_{k=1}^n \big( X_k - \mathbb{E} \big[ X_k \big]\big)\right| \geq t \Bigg) \leq 2\exp \Bigg[ - \frac{t^2}{2(\nu +ct)} \Bigg].
\end{equation*}
\end{theorem}

Denoting $G=n^{-1} (C_{11}^{\star})^{-1} X_1^T(\Lambda^\star)^{1/2}$, we get that
$$
\left(\frac1n (C_{11}^{\star})^{-1} X_1^T (\textbf{Y}-\boldsymbol{\lambda^{\star}})\right)_j=\sum_{k=1}^n G_{jk}\frac{(Y_k-\lambda_k^\star)}{\sqrt{\lambda_k^\star}}.
$$
Let us now apply the Bernstein's inequality to $X_k=G_{jk} Y_k/\sqrt{\lambda_k^\star}$ then
\begin{align*}
\sum_{k=1}^n\mathbb{E}[X_k^2]=\sum_{k=1}^n G_{jk}^2(1+\lambda_k^\star)& \leq n\sup_{1\leq k\leq n}(1+\lambda_k^\star)\;\|G\|_2^2 \\
&\leq n\sup_{1\leq k\leq n}(1+\lambda_k^\star)\;\rho(GG^T)\leq\sup_{1\leq k\leq n}(1+\lambda_k^\star)\;\rho((C_{11}^{\star})^{-1}).
\end{align*}
By Weyl's inequalities \cite[Corollary 4.3.15]{horn_johnson:2013} and \ref{th_2}, we get that with a probability tending to 1,
$$
\lambda_{min}(C_{11}^\star)\geq\lambda_{min}(C_{11})+\lambda_{min}(C_{11}^\star-C_{11})\geq\lambda_{min}(C_{11})-\rho(C_{11}^\star-C_{11})
\geq M_2-\|C_{11}^\star-C_{11}\|_2,
$$
Denoting $\bar{\lambda}=\max(\sup_{1\leq k\leq n}\lambda_k^\star,1)$ and $M_2'$ the positive constant such that $\rho((C_{11}^{\star})^{-1})\leq 1/M_2'$ for large enough $n$,
which exists by (\ref{eq:norme2_bound}), we get by \ref{th_2} that
$
\nu=2\bar{\lambda}/M_2'.
$
Observe that
\begin{equation}\label{eq:E_Xk_l}
  \sum_{k=1}^n \mathbb{E}\big[ |X_k|^{\ell} \big]=\sum_{k=1}^n \mathbb{E} \left[ \left| \frac{G_{jk}}{\sqrt{\lambda_k^{\star}}}Y_k \right|^{\ell} \right]
  = \sum_{k=1}^n \left( \frac{|G_{jk}|}{\sqrt{\lambda_k^{\star}}}\right)^{\ell} \sum_{i=1}^{\ell} {\lambda_k^{\star} }^{\ell} \bracenom{\ell}{i} ,
 \end{equation}
 where $\bracenom{\ell}{i}$ denotes the Stirling number of the second kind and the last equality is due to the definition of the $\ell$-th moment of a
 Poisson random variable. Then we have, for all $\ell \geq 3$, by Equation (2.3.8) of \cite{golub:96},
 \begin{align}\label{eq:Gjk}
   &  \sum_{k=1}^n \Bigg( \frac{|G_{jk}|}{\sqrt{\lambda_k^{\star}}}\Bigg)^{\ell} \sum_{i=1}^{\ell} {\lambda_k^{\star} }^{\ell} \bracenom{\ell}{i}  \leq \sum_{k=1}^n \frac{1}{n^{\ell/2}} \rho((C_{11}^\star)^{-1})^{\ell/2} \bar{\lambda}^{\ell/2} \sum_{i=1}^{\ell}  \bracenom{\ell}{i} \leq n^{1-\ell/2} M_2'^{-\ell/2} \bar{\lambda}^{\ell/2} \ell!
   \nonumber\\
 &= n^{1-\ell/2} M_2'^{1-\ell/2} \bar{\lambda}^{\ell/2-1}  \frac{\ell!}{2}\nu
   \leq\frac{\ell!}{2}\nu \left(\frac{\sqrt{\bar{\lambda}}}{\sqrt{n M_2'}}\right)^{\ell-2}=\frac{\ell!}{2} \nu c^{\ell-2},
 \end{align}
 with $c=\sqrt{\bar{\lambda}/(n M_2')}$ since $\sum_{i=1}^{\ell} \bracenom{\ell}{i} \leq \ell!$.

 Hence, for $t=n^{(c_1-1)/2}$,
 $$
 \frac{t^2}{2(\nu +ct)}=\frac{n^{c_1-1}}{2(\nu+n^{-1/2}\sqrt{\bar{\lambda}/M_2'}n^{(c_1-1)/2})}=O(n^{c_1/2}),
 $$
which gives the expected result.

 
Let us then prove that $\mathbb{P}(B_n^c)$ tends to zero as $n$ tends to infinity. 
By denoting
\begin{equation}
  \zeta = (\zeta_1, \ldots, \zeta_{p-q})^T = C_{21}C_{11}^{-1} W_1  -W_2 
  \quad \textrm{ and }\quad
  d = (d_1, \ldots, d_{p-q})^T =  C_{21}C_{11}^{-1} \text{sign}(\pmb{\beta}_1^{\star}),
  \label{def:zeta}
\end{equation}
we get that 
\begin{align}\label{eq:Bn_c_bound}
  \mathbb{P}(B_n^c) &= \mathbb{P} \left( |C_{21}C_{11}^{-1} W_1  -W_2| + |C_{21}C_{11}^{-1} R_1-R_2|+ \frac{\alpha}{2n} \left|  C_{21}C_{11}^{-1} \text{sign}(\pmb{\beta}_1^{\star})\right|
                     >  \frac{\alpha}{2n}   \right) \nonumber\\
  & \leq \sum_{j=1}^{p-q} \mathbb{P}\left(|\zeta_j| + |(C_{21}C_{11}^{-1} R_1-R_2)_j|+\frac{\alpha}{2n} \big| d_j \big|\geq \frac{\alpha}{2n}  \right)\nonumber\\
  &  \leq \sum_{j=1}^{p-q} \left\{\mathbb{P}\left(|\zeta_j|\geq  \frac{\alpha}{6n}\right)+\mathbb{P}\left(|(C_{21}C_{11}^{-1} R_1-R_2)_j|\geq \frac{\alpha}{6n}\right)
    +\mathbb{P}\left(\big| d_j \big|\geq  \frac{1}{3}\right)\right\}.
  \end{align}
  By the strong irrepresentable condition \eqref{eq:IC}, we get that 
  $\sum_{j=1}^{p-q}\mathbb{P}\left(\big| d_j \big|\geq  \frac{1}{3}\right)=o(1)$.
Let us now prove that the second term in the r.h.s of (\ref{eq:Bn_c_bound}) tends to 0 as $n$ tends to infinity.
  Observing that $R_1$ defined in (\ref{eq:R1_2}) satisfies: $R_1= C_{11}\left( \pmb{\beta}^\star_1 - \pmb{\Tilde{\beta}}_1  \right) + C_{12}\left( \pmb{\beta}^\star_2 - \pmb{\Tilde{\beta}}_2  \right)$, we get by the Cauchy-Schwarz inequality that
 \begin{align*}
\left | (C_{21}C_{11}^{-1} R_1)_j \right | &\leq \sqrt{p-q} \left \| C_{21}C_{11}^{-1} R_1\right \|_2 \\
&\leq \sqrt{p-q} \left \|C_{21}( \pmb{\beta}^\star_1 - \pmb{\Tilde{\beta}}_1 ) \right\|_2 + \sqrt{p-q} \left \| C_{21}C_{11}^{-1} C_{12} \left(\pmb{\beta}^\star_2 - \pmb{\Tilde{\beta}}_2  \right) \right \|_2.
\end{align*}
By \ref{th_5}, \ref{th_2} and \ref{th_3}, $\left | (C_{21}C_{11}^{-1} R_1)_j \right | =O_P(n^{-1})$ for all $j$. Using similar arguments  we get $\left | (R_2)_j \right | =O_P(n^{-1})$.  Then since $\alpha=O\left(n^{(c_2+1)/2}\right)$, 
we get that the second term in the r.h.s of (\ref{eq:Bn_c_bound}) tends to 0 as $n$ tends to infinity. 

Let us now prove that the first term in the r.h.s of (\ref{eq:Bn_c_bound}) tends to 0 as $n$ tends to infinity. By \eqref{def:zeta}, \eqref{eq:def_C_W}, \eqref{eq_eps} and (\ref{eq:Y_X}), $\zeta$ can be rewritten as follows for all $j\in\{1,\ldots,p-q\}$
\begin{equation}\label{eq:zeta}
\zeta_j = (C_{21}C_{11}^{-1} W_1-W_2)_j =  \left(C_{21}\xi- \frac{\mathcal{X}_2^T}{n} \Tilde{\varepsilon}\right)_j=\left(C_{21}\xi-\frac{X_2^T}{n} (\textbf{Y}-\boldsymbol{\lambda^{\star}})\right)_j,
\end{equation}
where we recall that $\boldsymbol{\lambda^{\star}}$ denotes the column vector of the $\lambda_i^\star$'s and $\xi = C_{11}^{-1} W_1$.

Let us consider the first term in \eqref{eq:zeta} and prove that $\mathbb{P} \Bigg( |(C_{21}\xi)_j| \geq\frac{\alpha}{12n}\Bigg)$ tends to $0$ as $n$ tends to infinity. Let us note that the term $\xi_j$ is handled previously in the proof concerning $A_n^c$. Therefore we use the same arguments and adapt it to the term 
\begin{equation}\label{eq:c21xi}
(C_{21}\xi)_j= \left(\frac1n C_{21}C_{11}^{-1} \mathcal{X}_1^T \boldsymbol{\Tilde{\Lambda}}^{-1/2}(\textbf{Y}-\boldsymbol{\lambda^{\star}})\right)_j
  +\left(\frac1n C_{21}C_{11}^{-1} \mathcal{X}_1^T \boldsymbol{\Tilde{\Lambda}}^{-1/2}(\boldsymbol{\lambda^{\star}}-\boldsymbol{\tilde{\lambda}})\right)_j.
  \end{equation}
The second term of \eqref{eq:c21xi} is still $O_P(1/\sqrt{n})$ since in \eqref{eq:rho}, $\rho(C_{21}C_{11}^{-1}C_{21}^T)=\rho(C_{21}C_{11}^{-1}C_{12})$ is bounded by using \ref{th_2} and \ref{th_3}. The first term of \eqref{eq:c21xi} can be decomposed in the same way as \eqref{eq:2eme_terme}:
\begin{align}\label{eq:2eme_terme:c12xi}
&\left(\frac1n C_{21}C_{11}^{-1} \mathcal{X}_1^T \boldsymbol{\Tilde{\Lambda}}^{-1/2}(\textbf{Y}-\boldsymbol{\lambda^{\star}})\right)_j\nonumber\\
&=\left(\frac1n C_{21}^\star(C_{11}^{\star})^{-1} X_1^T (\textbf{Y}-\boldsymbol{\lambda^{\star}})\right)_j+
\left(\frac1n (C_{21}C_{11}^{-1} -C_{21}^\star(C_{11}^{\star})^{-1}) X_1^T (\textbf{Y}-\boldsymbol{\lambda^{\star}})\right)_j.
\end{align}
The second term of \eqref{eq:2eme_terme:c12xi} is bounded in the same manner as the second term of \eqref{eq:2eme_terme} was handled. 
By observing that
  \begin{align*}
C_{21}C_{11}^{-1} -C_{21}^\star(C_{11}^{\star})^{-1}&=(C_{21}-C_{21}^\star)C_{11}^{-1} + C_{21}^\star(C_{11}^{-1} -(C_{11}^\star)^{-1} )\\
&=(C_{21}-C_{21}^\star)C_{11}^{-1} + (C_{21}^\star-C_{21})(C_{11}^{-1} -(C_{11}^\star)^{-1} )+C_{21}(C_{11}^{-1} -(C_{11}^\star)^{-1} ),
\end{align*}
and by using (\ref{eq:norme2_bound_C11-1}), \ref{th_1}, \ref{th_2}, \ref{th_3} and the fact that
\begin{equation}\label{eq:C21}
  \|C_{21}^\star-C_{21}\|_2=O_P(1/n^2),
\end{equation}
where we used the same arguments as in (\ref{eq:norme2_bound}) we get that
  $\|C_{21}C_{11}^{-1} -C_{21}^\star(C_{11}^{\star})^{-1}\|_2=O_P(1/n^2)$. Thus, the second term of\eqref{eq:2eme_terme:c12xi} is $O_P(1/n^2)$.
To handle the first term of \eqref{eq:2eme_terme:c12xi}, let us apply the Bernstein's inequality to $X_k=\widetilde{G}_{jk} Y_k/\sqrt{\lambda_k^\star}$ with $\widetilde{G}=n^{-1}C_{21}^\star(C_{11}^{\star})^{-1} X_1^T(\Lambda^\star)^{1/2}$. Then
$$
\sum_{k=1}^n\mathbb{E}[X_k^2]=\sum_{k=1}^n \widetilde{G}_{jk}^2(1+\lambda_k^\star)\leq n\sup_{1\leq k\leq n}(1+\lambda_k^\star)\;\rho(\widetilde{G}\widetilde{G}^T)\leq\sup_{1\leq k\leq n}(1+\lambda_k^\star)\;\rho(C_{21}^{\star}(C_{11}^\star)^{-1}(C_{21}^{\star})^T).
$$
By Weyl's inequalities \cite[Corollary 4.3.15]{horn_johnson:2013}
  \begin{align*}
    & \lambda_{max}(C_{21}^{\star}(C_{11}^\star)^{-1}(C_{21}^{\star})^T)\leq
      \lambda_{max}((C_{21}^{\star}-C_{21})((C_{11}^{\star})^{-1}-C_{11}^{-1})((C_{21}^{\star})^T-C_{21}^T))\\
&  +\lambda_{max}((C_{21}^{\star}-C_{21})((C_{11}^{\star})^{-1}-C_{11}^{-1})C_{21}^T)+\lambda_{max}((C_{21}^{\star}-C_{21})C_{11}^{-1}((C_{21}^{\star})^T-C_{21}^T))\\
&   +\lambda_{max}((C_{21}^{\star}-C_{21})C_{11}^{-1}C_{21}^T)+\lambda_{max}(C_{21}((C_{11}^{\star})^{-1}-C_{11}^{-1})((C_{21}^{\star})^T-C_{21}^T))\\
&  +\lambda_{max}(C_{21}((C_{11}^{\star})^{-1}-C_{11}^{-1})C_{21}^T)+\lambda_{max}(C_{21}C_{11}^{-1}((C_{21}^{\star})^T-C_{21}^T))
  +\lambda_{max}(C_{21}C_{11}^{-1}C_{21}^T).
  \end{align*}

By (\ref{eq:norme2_bound_C11-1}), (\ref{eq:C21}), \ref{th_2} and \ref{th_3}, we get that for a large enough $n$,
  $\rho(C_{21}^{\star}(C_{11}^\star)^{-1}(C_{21}^{\star})^T)\leq M_3'$, where $M_3'$ is a positive constant. Hence, $\nu=(2\bar{\lambda})/M_3'$. Using the same
  bounds as in (\ref{eq:E_Xk_l}) and (\ref{eq:Gjk}), we get that $c=\sqrt{\bar{\lambda}/(nM'_3)}$. Thus, with $t=n^{(c_2-1)/2}$ in the Bernstein's inequality
  $$
  \frac{t^2}{2(\nu +ct)}=\frac{n^{c_2-1}}{2(\nu+n^{-1/2}\sqrt{\bar{\lambda}/M_3'}n^{(c_2-1)/2})}=O(n^{c_2/2}).
$$
Consequently, we conclude that for $\alpha=O\left(n^{(c_2+1)/2}\right)$, $\mathbb{P} \Bigg( |(C_{21}\xi)_j| \geq\frac{\alpha}{12n}\Bigg)$ tends to $0$ as $n$ tends to infinity.

Finally, let us consider the second term in \eqref{eq:zeta} and prove that $\mathbb{P} \Bigg( \left|\left(\frac{X_2^T}{n} (\textbf{Y}-\boldsymbol{\lambda^{\star}})\right)_j\right| \geq\frac{\alpha}{12n}\Bigg)$ tends to $0$ as $n$ tends to infinity.
Denoting $H=n^{-1} X_2^T(\Lambda^\star)^{1/2}$, we observe that
$$
\left(\frac1n X_2^T (\textbf{Y}-\boldsymbol{\lambda^{\star}})\right)_j=\sum_{k=1}^n H_{jk}\frac{(Y_k-\lambda_k^\star)}{\sqrt{\lambda_k^\star}}
$$
and we apply the Bernstein's inequality to $X_k=H_{jk} Y_k/\sqrt{\lambda_k^\star}$. Then we have
\begin{align*}
\sum_{k=1}^n\mathbb{E}[X_k^2]=\sum_{k=1}^n H_{jk}^2(1+\lambda_k^\star)& \leq n\sup_{1\leq k\leq n}(1+\lambda_k^\star)\;\|H\|_2^2 \\
&\leq n\sup_{1\leq k\leq n}(1+\lambda_k^\star)\;\rho(HH^T)\leq\sup_{1\leq k\leq n}(1+\lambda_k^\star)\;\rho(C_{22}^{\star}).
\end{align*}
By Weyl's inequalities \cite[Corollary 4.3.15]{horn_johnson:2013}
$$
\lambda_{max}(C_{22}^{\star})\leq\lambda_{max}(C_{22})+\lambda_{max}(C_{22}^{\star}-C_{22}).
$$
By using the same arguments as in (\ref{eq:norme2_bound}) and \ref{th_3}, we get that
$\rho(C_{22}^{\star})\leq M_4'$ for large enough $n$ and by denoting $\bar{\lambda}=\max(1,\sup_{1\leq k\leq n}\lambda_k^\star)$, we get  that
$\nu=2\bar{\lambda}M_4'$. Then, using exactly the same argument as before we get
$$ \sum_{k=1}^n \mathbb{E}\big[ |X_k|^{\ell} \big]\leq \frac{\ell!}{2} \nu c^{\ell-2},$$ 
with $c=\sqrt{\bar{\lambda}/(n M_4')}$. Hence, for $t=n^{(c_2-1)/2}$ in the Bernstein's inequality,
  we obtain the expected result, which concludes the proof.
 



\bibliographystyle{elsarticle-num} 
\bibliography{biblio}

\begin{thebibliography}{1}
\expandafter\ifx\csname url\endcsname\relax
  \def\url#1{\texttt{#1}}\fi
\expandafter\ifx\csname urlprefix\endcsname\relax\def\urlprefix{URL }\fi
\expandafter\ifx\csname href\endcsname\relax
  \def\href#1#2{#2} \def\path#1{#1}\fi

\bibitem{friedman_hastie_tibshirani:2010}
J.~H. Friedman, T.~Hastie, R.~Tibshirani, Regularization paths for generalized
  linear models via coordinate descent, Journal of Statistical Software 33~(1)
  (2010).

\bibitem{wu2017diversity}
H.~Wu, L.~Yang, L.-L. Chen, The diversity of long noncoding {RNA}s and their
  generation, Trends in genetics 33~(8) (2017) 540--552.

\bibitem{tibshirani:1996}
R.~Tibshirani, Regression shrinkage and selection via the lasso, Journal of the
  Royal Statistical Society: Series B (Methodological) 58~(1) (1996) 267--288.

\bibitem{zhao:2006}
P.~Zhao, B.~Yu, On model selection consistency of lasso, The Journal of Machine
  Learning Research 7 (2006) 2541--2563.

\bibitem{giraud:2021}
C.~Giraud, Introduction to high-dimensional statistics, Chapman and Hall/CRC,
  2021.

\bibitem{golub:96}
H.~Golub, G.\, F.~Van~Loan, C.\, Matrix Computations, JHU Press, 1996.

\bibitem{horn_johnson:2013}
R.~A. Horn, C.~R. Johnson, Matrix Analysis, Cambridge University Press, 2013.

\bibitem{boucheron_lugosi_massart:2013}
S.~Boucheron, G.~Lugosi, P.~Massart, Concentration inequalities: A
  nonasymptotic theory of independence, Oxford university press, 2013.

\end{thebibliography}






\end{document}